\documentclass[twocolumn]{autart}    

\usepackage{graphicx}          

\usepackage{colortbl}

\newtheorem{theorem}{Theorem}
\newtheorem{corollary}{Corollary}

\newtheorem{lemma}{Lemma}

\newtheorem{definition}{Definition}

\begin{document}

\begin{frontmatter}

\title{Exponential Stabilization for It\^{o} Stochastic Systems with Multiple Input Delays
\thanksref{footnoteinfo}} 

\thanks[footnoteinfo]{This work is supported by the National Natural Science Foundation of China (61633014, 61573221) and the
Qilu Youth Scholar Discipline Construction Funding from Shandong University.
Corresponding author H. Zhang.}

\author[sdu]{Juanjuan Xu}\ead{jnxujuanjuan@163.com},    
\author[sdu]{Huanshui Zhang}\ead{hszhang@sdu.edu.cn}  

\address[sdu]{School of Control Science and Engineering, Shandong University, Jinan, Shandong, P.R. China 250061}  


\begin{abstract}                          
In this paper, we study the stabilization problem for the It\^{o} systems with both multiplicative noise and multiple delays
which exist widely in applications such as networked control systems. Sufficient and necessary conditions are obtained for the
exponential stabilization problem of It\^{o} stochastic systems with multiple delays. On one hand, we derive the solvability of the modified
Riccati equation in case of the mean-square exponential stabilization. On the other hand, the mean-square exponential stabilization is
guaranteed by the solvability of a modified Riccati equation. A novel stabilizing controller is shown in the feedback from of the conditional
expectation in terms of the modified algebraic Riccati equation.
The main technique is to reduce the original system with multiple delays to a pseudo delay-free system.
\end{abstract}

\end{frontmatter}

\section{Introduction}
The mathematical models described by delayed differential equations are ubiquitous and have wide applications
in physics, engineering, communication, biology and so on \cite{d1}. As is known, time delays usually degrade the system performance,
and are the source of instability, and even lead to the occurrence of chaos phenomenon. So study on the stabilization problem of time-delay
system is of great significance. Some essential progress has been made on the optimal control and
stabilization problems for time delay systems, see \cite{d2}, \cite{smith} and references therein.
In particular, \cite{smith} designs a predictor-like controller which reduces the original delayed system to delay-free one.
By virtue of the predictor-based technique, the problem for systems with more general delays has been studied in \cite{Artstein}-\cite{Olbrot}.
The linear quadratic regulation (LQR) problem for systems with multiple input delays was solved in \cite{zhang2} by establishing a duality between the LQR problem
and a smoothing problem. The optimal controller is presented using a Riccati equation.
\cite{tadmor1}-{\cite{tadmor2}} studied the $H_\infty$ preview control problem and presented the necessary and sufficient
solvability conditions in terms of a standard algebraic Riccati equation and a nonstandard $H_\infty$-like algebraic Riccati equation.
The aforementioned results are only related to the deterministic system and more details are referred to the
survey paper \cite{d2}.

Considering the accuracy requirement to the system in applications, it is necessary to take the uncertainty into consideration.
One of the most popular models is the stochastic differential equation motivated by Brownian motion.
When the stochastic system is delay-free, \cite{rami} presents some sufficient and necessary conditions
for the mean-square stabilization. There have also been many important developments when both delay and uncertainty are considered, especially the noise is multiplicative, e.g.,
\cite{s1}, \cite{s3}, \cite{zdwang} and references therein. Noting that most results in the literature depend on the linear matrix inequality (LMI) to characterize the sufficient conditions for the stabilization.
For instance, \cite{zdwang} investigated the stochastic stabilization problem for a class of bilinear continuous time-delay uncertain systems
with Markovian jumping parameters. Sufficient conditions were established to guarantee the existence of desired
robust controllers, which are given in terms of the solutions to a set of LMIs, or coupled quadratic matrix
inequalities. \cite{xie} considered a class of large-scale interconnected bilinear stochastic systems with time delays and time-varying parameter
uncertainties and robust stability analysis was given in terms of a set of LMIs.
In addition, some convergence theorems have been given in the literature. For example, \cite{mao1}-\cite{mao2} investigated the LaSalle-type asymptotic
convergence theorems for the solutions of stochastic differential delay equations. More recently, some substantial progress for the optimal LQ control has been made by proposing the approach
of solving the forward and backward differential/difference equations (FBDEs). See \cite{hszhang} and \cite{hszhang1} for details. 
However, the stabilization problem for It\^{o} stochastic systems with multiple delays have not yet been completely solved. The main obstacles are that the problem is in fact infinite dimensional and the classical controller such as
current feedback form only leads to sufficient conditions which may be delay-dependent. 

Inspired by the work \cite{hszhang1}, we shall study the stochastic system with multiple delays. The main contribution is two-fold.
Firstly, we derive the solvability of the modified Riccati equation in case of the mean-square exponential stabilization.
Secondly, we obtain that the mean-square exponential stabilization can be
guaranteed by the solvability of a modified Riccati equation.
A novel stabilizing controller is shown in the feedback from of the conditional
expectation in terms of the modified algebraic Riccati equation.
The main technique is to reduce the original system with multiple delays to a pseudo delay-free system.

The rest of the paper is formulated as follows: Section 2 illustrates the studied problem.
The system is reduced to a pseudo delay-free system and the optimization problems of the reduced system are studied in Section 3.
Sufficient and necessary conditions are given in Section 4 for the exponential mean-square stabilization of the system.
Some concluding remarks are shown in the last section.

\textbf{Notation.} $R^n$
denotes the family of $n$-dimensional vectors; $x'$ denotes the
transpose of $x$; and a symmetric matrix $M>0\ (\geq 0)$ is
strictly positive-definite (positive semi-definite). $(\Omega, \mathcal{F}, \mathcal{P}, \mathcal{F}_t|_{t\geq 0})$ is a complete
stochastic basis so that $\mathcal{F}_0$ contains all P-null
elements of $\mathcal{F},$ and the filtration is generated by the
standard Brownian motion $\{w(t)\}_{t\geq 0}.$ $\hat{x}(t|s)\doteq
E[x(t)|\mathcal{F}_{s}]$ denotes the conditional expectation with
respect to the filtration $\mathcal{F}_{s}.$ We simply denote
$E_t(\cdot)=E[\cdot|\mathcal{F}_{t}],$ and $\langle\cdot,\cdot\rangle$ denotes the
inner product in Hilbert space. The following sets are useful
throughout the paper:
\begin{eqnarray}
\bar{C}_{[-h,0)}&=&\{\varphi(t):
[-h,0)\rightarrow R^m \mbox{~is~continuous~and~}
\nonumber\\
&&\sup_{-h\leq
t<0}\|\varphi(t)\|<\infty\},\nonumber\\
L_{\mathcal{F}}^2(0,T;R^m)&=&\{\varphi(t)_{t\in[0,T]}
\mbox{~is~an~}\mathcal{F}_t-\mbox{adapted~stochastic}\nonumber\\
&&\mbox{~process}~s.t.~
E\int_0^T\|\varphi(t)\|^2dt<\infty\}.\nonumber
\end{eqnarray}

\section{Problem Formulation}

Consider the It\^{o} stochastic systems with multiple input delays:
\begin{eqnarray}
dx(t)&=&\Big(Ax(t)+\sum_{i=0}^{r}B_iu(t-h_i)\Big)dt\nonumber\\
&&+\sum_{i=0}^{r}\bar{B}_iu(t-h_i)dw_i(t),\label{i1}
\end{eqnarray}
where $x(t)\in R^n$ is the state, $u(t)\in R^m$ is the control
input, $h_0=0, h_i>0,i=1,\cdots,r$ represent the input delays. $w_i(t),i=1,\cdots,r$ is independent one-dimension
standard Brownian motion. $A,B_i,\bar{B}_i$ are constant matrices with
compatible dimensions.  The initial conditions are chosen as
$x(0)=x_0$ and $u(\tau)=\mu(\tau)\in \bar{C}_{[-h_r,0)}.$

\emph{Remark 1.} The system (\ref{i1}) has wide applications in network control systems. In particular, consider the continuous-time LTI system
with both random input gains and multiple input delays as shown in Fig. \ref{fig1}:
\begin{figure}
\begin{center}
\includegraphics[height=4cm]{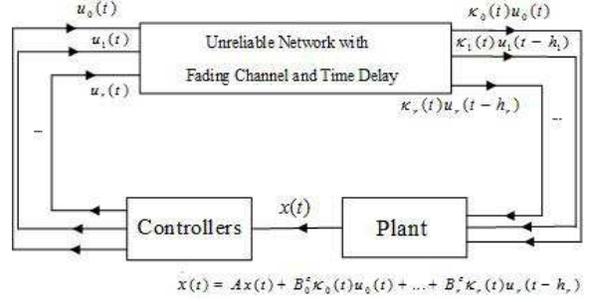}  
\caption{Continuous-time LTI system
with both random input gains and multiple input delays}  
\label{fig1}                                 
\end{center}                                 
\end{figure}

\begin{eqnarray}
\dot{x}(t)&=&Ax(t)+B^c_0\kappa_0(t)u_0(t)+B^c_1\kappa_1(t)u_1(t-h_1)+\cdots\nonumber\\
&&+B^c_r\kappa_r(t)u_r(t-h_r).\label{ex1}
\end{eqnarray}
where $x(t)\in R^n$ is the state, $u_i(t)\in R^m$ is the $i$th control
input, $h_i>0,i=1,\cdots,r$ represent the input delays.
$\kappa_i(t)=\mu_i+\xi_i(t)$ where $\mu_i$ is a real positive constant and $\xi_i(t)$ is a zero-mean white noise with autocorrelation
$E[\xi_i(t)\xi_i(t+\tau)] = \sigma_i^2\delta(\tau)$.
By denoting $u(t)=\left[
               \begin{array}{ccc}
                 u_0(t) & \cdots & u_r(t) \\
               \end{array}
             \right]'$
and $B_i=\left[
             \begin{array}{ccccccc}
               0 & \cdots & 0 & B_i^c & 0 & \cdots & 0 \\
             \end{array}
           \right]$ for $i=0,\cdots,r,$ (\ref{ex1}) can be rewritten as
\begin{eqnarray}
\dot{x}(t)&=&Ax(t)+B_0\kappa_0(t)u(t)+B_1\kappa_1(t)u(t-h_1)+\cdots\nonumber\\
&&+B_r\kappa_r(t)u(t-h_r).\label{ex2}
\end{eqnarray}
(\ref{ex2}) can then be reformulated as a standard It\^{o} form by using $\kappa_i(t)=\mu_i+\xi_i(t)$:
\begin{eqnarray}
dx(t)&=&Ax(t)dt+B_0u(t)[\mu_0dt+\sigma_0dw_0(t)]\nonumber\\
&&+B_1u(t-h_1)[\mu_1dt+\sigma_1dw_1(t)]\nonumber\\
&&+\cdots+B_ru(t-h_r)[\mu_rdt+\sigma_rdw_r(t)]\nonumber\\
&=&\Big(Ax(t)+\sum_{i=0}^{r}\mu_iB_iu(t-h_i)\Big)dt\nonumber\\
&&+\sum_{i=0}^{r}\sigma_i\bar{B}_iu(t-h_i)dw_i(t).\nonumber
\end{eqnarray}
This is a special case of systems (\ref{i1}).

We now define the stabilization and exponential stabilization for system (\ref{i1}).
\begin{definition}\label{d1}
System (\ref{i1}) is mean-square stabilizable if there
exists an $\mathcal{F}_t$-adapted controller $u(t)$ in the form of
\begin{eqnarray}
Lx(t)+\int_t^{t+h_r}L(s)u(s-h_r)ds\label{i2}
\end{eqnarray}
where $L$ is a constant matrix and $L(s)$ is a time-varying matrix with compatible dimensions
such that the closed-loop system satisfies
\begin{eqnarray}
\lim_{t\rightarrow\infty}E\|x(t)\|^2=0 ~~\mbox{and}~~ \lim_{t\rightarrow\infty}E\|u(t)\|^2=0\nonumber
\end{eqnarray}
for any $x_0$ and any $\mathcal{F}_t$-adapted controller
$u(t),t\leq h_r.$
\end{definition}

\begin{definition}\label{d2}
System (\ref{i1}) is mean-square exponentially stabilizable if there
exists an $\mathcal{F}_t$-adapted controller $u(t)$ in the form of
(\ref{i2}) and a positive constant $\alpha$
such that the closed-loop system satisfies
\begin{eqnarray}
\lim_{t\rightarrow\infty}e^{\alpha t}E\|x(t)\|^2=0~~ \mbox{and}~~\lim_{t\rightarrow\infty}e^{\alpha t}E\|u(t)\|^2=0\nonumber
\end{eqnarray}
for any $x_0$ and any $\mathcal{F}_t$-adapted controller
$u(t),t\leq h_r.$
\end{definition}

The aim of this paper is stated as follows.
\\
{\em Problem }: Find the sufficient and necessary conditions for system (\ref{i1}) to be exponentially stabilized by a controller
in the form of (\ref{i2}) following Definition \ref{d2}.

The outline of the solvability to \emph{Problem} is as follows: Firstly, we convert the original stochastic system with multiple
input delays into a pseudo delay-free system where the delays are involved in the Brownian motions rather than the control input.
Secondly, we solve finite-horizon optimization problems with a standard cost function and a discounted cost function
subject to the pseudo delay-free system in terms of modified differential Riccati equations. Finally, the sufficient and necessary
conditions for the exponential stabilization are characterized by the corresponding modified algebraic Riccati equation.

\section{Reduction of the original system into a pseudo delay-free system}

We firstly transform the original system (\ref{i1}) into a pseudo delay-free system. To this end, we define
\begin{eqnarray}
y(t)&=&x(t)+\sum_{i=1}^{r}\int_t^{t+h_i}e^{A(t-s)}B_iu(s-h_i)ds\nonumber\\
&&+\sum_{i=1}^{r}\int_t^{t+h_i}e^{A(t-s)}\bar{B}_iu(s-h_i)dw_i(s).\label{i3}
\end{eqnarray}
\begin{lemma}
$y(t)$ defined by (\ref{i3}) satisfies the dynamic
\begin{eqnarray}
dy(t)
&=&\Big(Ay(t)+\sum_{i=0}^{r}e^{-Ah_i}B_iu(t)\Big)dt\nonumber\\
&&+\sum_{i=0}^{r}e^{-Ah_i}\bar{B}_iu(t)dw_i(t+h_i).\label{i4}
\end{eqnarray}
\end{lemma}
\emph{Proof.}
By taking It\^{o}'s formula to $y(t)$ and using (\ref{i1}), it is obtained that
\begin{eqnarray}
dy(t)&=&\Big(Ax(t)+\sum_{i=0}^{r}B_iu(t-h_i)\Big)dt\nonumber\\
&&+\sum_{i=0}^{r}\bar{B}_iu(t-h_i)dw_i(t)+\sum_{i=1}^{r}e^{-Ah_i}B_iu(t)dt\nonumber\\
&&-\sum_{i=1}^{r}B_iu(t-h_i)dt+\sum_{i=1}^{r}e^{-Ah_i}\bar{B}_iu(t)dw_i(t+h_i)\nonumber\\
&&-\sum_{i=1}^{r}\bar{B}_iu(t-h_i)dw_i(t)\nonumber\\
&&+A\Big(\sum_{i=1}^{r}\int_t^{t+h_i}e^{A(t-s)}B_iu(s-h_i)ds\nonumber\\
&&+\sum_{i=1}^{r}\int_t^{t+h_i}e^{A(t-s)}\bar{B}_1u(s-h_i)dw_i(s)\Big)dt\nonumber\\
&=&\Big(Ay(t)+\sum_{i=0}^{r}e^{-Ah_i}B_iu(t)\Big)dt\nonumber\\
&&+\sum_{i=0}^{r}e^{-Ah_i}\bar{B}_iu(t)dw_i(t+h_i).\nonumber
\end{eqnarray}
This completes the proof.

\emph{Remark 2.} Noting that there exists no delay in the control input $u(t).$ However,
the delays $h_i,i=1,\cdots,r$ are involved in the Brownian motions $w_i.$ Thus, we call the system (\ref{i4}) as a pseudo delay-free system.

Define a new $\sigma$-algebraic $\mathcal{G}_t=\{w_i(s+h_i),i=0,1,\ldots,r, s\leq t\}.$
Then it holds that $\mathcal{F}_t\subseteq \mathcal{G}_t\subseteq \mathcal{F}_{t+h_r}.$
From (\ref{i4}), we have $y(t)$ is $\mathcal{G}_t$-adapted. In addition, considering Definition \ref{d1} and (\ref{d2}),
the controller $u(t)$ is $\mathcal{F}_t$-adapted. For convenience of the future use, it is simply denoted that $B=\sum_{i=0}^{r}e^{-Ah_i}B_i.$

\subsection{Finite-horizon optimal control problem of pseudo delay-free system}

We then study the finite-horizon optimization problem of minimizing the standard linear quadratic cost function subject to (\ref{i4}):
\begin{eqnarray}
J_T&=&E\Big\{
\int_{0}^{T }\Big(y'(t)Qy(t)+u(t)'Ru(t)\Big)dt\nonumber\\
&&+y'(T)Hy(T)\Big\},\label{2}
\end{eqnarray}
where $H$ is semi-positive definite matrix of compatible dimension.

Noting that the new state $y(t)$ is $\mathcal{G}_t$-adapted rather than $\mathcal{F}_t$-adapted, we define the admissible control set as
\begin{eqnarray}
\mathcal{U}_{ad}&=&\lbrace u(t)\in L_{\mathcal{F}}^2(0,\infty;R^m):u(t)=M(t)\hat{y}(t|t)\},
\label{i7}
\end{eqnarray}
where $M(t)$ is time-varying matrices with compatible dimension and
\begin{eqnarray}
\hat{y}(t|t)&=&E[y(t)|\mathcal{F}_t]\nonumber\\
&=&x(t)+\sum_{i=1}^{r}\int_t^{t+h_i}e^{A(t-s)}B_iu(s-h_i)ds.\nonumber
\end{eqnarray}

Following \cite{zhang3}, the stochastic maximum principle can be immediately obtained.
\begin{lemma} \label{mp}
The optimal solution to minimize (\ref{2}) subject to (\ref{i4}) satisfies
\begin{eqnarray}
0&=&Ru(t)+E[B'p(t)+\sum_{i=0}^{r}\bar{B}_i'e^{-A'h_i}q_i(t)|\mathcal{F}_{t}],\label{4}
\end{eqnarray}
where $(p(t),q(t))$ is the solution of the backward stochastic
differential equation (BSDE):
\begin{eqnarray}\label{5}
\left\{
  \begin{array}{ll}
    dp(t)=-[A'p(t)+Qy(t)]dt+\sum_{i=0}^{r}q_i(t)dw_i(t), & \hbox{} \\
    p(T)=Hx(T). & \hbox{}
  \end{array}
\right.\hspace{-0.8cm}
\end{eqnarray}
while $y(t)$ obeys (\ref{i4}) and $H$ is defined in (\ref{2}).
\end{lemma}

Based on Lemma \ref{mp}, the explicit solvability of forward and backward stochastic differential equations (\ref{i4}), (\ref{4}) and (\ref{5})
is the key to the derivation of the optimal solution. To this end, we define the modified differential Riccati equation:
\begin{eqnarray}
-\frac{d}{dt}\hat{P}(t)&=&\hat{P}(t)A+A'\hat{P}(t)+Q-\Pi(t,t),\label{x1}
\end{eqnarray}
and
\begin{eqnarray}
P(t)
&=&\hat{P}(t)+\int_0^{h_r}e^{A'\theta}\Pi(t+\theta,t+\theta)e^{A\theta}d\theta
,\label{x2}
\end{eqnarray}
where
\begin{eqnarray}
\Pi(t,t)&=&K'(t) \Omega(t) K(t), \label{x3}\\
\Omega(t) &=&R+\sum_{i=0}^{r}\bar{B}_i'e^{-A'h_i}P(t)e^{-Ah_i}\bar{B}_i,\label{x4}\\
K(t)&=&-\Omega^{-1}(t)B'\hat{P}(t),\label{x5}
\end{eqnarray}
with the terminal values $\hat{P}(T)=P(T)=H$ for $H$ defined in (\ref{2}).

\begin{lemma}\label{lem2}
The equation (\ref{x1})-(\ref{x5}) is equivalent to the following equations:
\begin{eqnarray}
-\dot{P}(t)&=&P(t)A+A'P(t)+Q-e^{A'h_r}\Pi(t+h_r,t+h_r)\nonumber\\
&&\times e^{Ah_r},
\label{6}
\end{eqnarray}
while $\Pi(t+h_r,t+h_r)$ is given by
\begin{eqnarray}
\Pi(t,t)&=&K'(t)\Omega(t)K(t),\label{8}\\
\Omega(t)&=&R+\sum_{i=0}^{r}\bar{B}_i'e^{-A'h_i}P(t)e^{-Ah_i}\bar{B}_i,\label{e140}\\
K(t)&=&-\Omega^{-1}(t)\Big[B'P(t)
-B'\int_{0}^{h_r}e^{A'\theta}\Pi(t+\theta,t+\theta)\nonumber\\
&&\times e^{A\theta}d\theta\Big],\label{e15}
\end{eqnarray}
with terminal values $P(T)=H$ and $\Pi(T,T+\theta)=0$ for $\theta\in (0,h_r].$
\end{lemma}

Proof. The equivalence can be established by similar discussions to Remark 5 in \cite{hszhang1}. So we omit it.

We now present the optimal solution of the finite-horizon linear quadratic optimal control problem by using the solution to (\ref{x1})-(\ref{x5}).
\begin{lemma}\label{lem1}
Assume that the modified Riccati equation (\ref{x1})-(\ref{x5}) admits a solution such that the matrix $\Omega(t)>0$, then there exists
a unique solution to the problem of minimizing (\ref{2}) subject to the system (\ref{i4}) and the optimal controller is given by
\begin{eqnarray}
u(t)&=&K(t)\hat{y}(t|t).\label{contr}
\end{eqnarray}
The optimal cost is as
\begin{eqnarray}
J_T^{\ast}&=&E\Big(y'(0)P(0)y(0)-y'(0)\int_{0}^{h_r}
\Pi(0,\theta)\hat{y}(0|\theta)d\theta\Big).\label{zz18}
\end{eqnarray}
\end{lemma}

Proof. The proof is presented in Appendix \ref{ap2}.

As a byproduct of Lemma \ref{lem1} which is useful in the stabilization, we further state the following results.

\begin{corollary}\label{cor2}
Under the same conditions in Lemma \ref{lem1} and let the controller satisfy that $u(t)=0$ for $t\in[-h_r,0).$
Then there exists a unique solution to the problem of minimizing (\ref{2}) subject to the system (\ref{i4}). The optimal controller is given by
(\ref{contr}) for $t\geq 0$ and the optimal cost is as
\begin{eqnarray}
J_T^{\ast}&=&E\Big(x_0'\hat{P}(0)x_0\Big).\label{d10}
\end{eqnarray}
\end{corollary}

Proof. Since $u(t)=0$ for $t\in[-h_r,0),$ then $y(0)=x(0).$ Thus the optimal cost becomes $J_T^{\ast}=E\Big(x_0'\hat{P}(0)x_0\Big)$ from (\ref{zz18}).

\begin{corollary}\label{cor1}
Under the same conditions in Lemma \ref{lem1} and let the controller satisfy that $u(t)=0$ for $t\in[-h_r,h_r).$
Then there exists a unique solution to the problem of minimizing (\ref{2}) subject to the system (\ref{i4}). The optimal controller is given by
(\ref{contr}) for $t\geq h_r$ and the optimal cost is as
\begin{eqnarray}
J_T^{\ast}&=&E\Big(x_0'P(0)x_0\Big).\label{i9}
\end{eqnarray}
\end{corollary}

Proof. Since $u(t)=0$ for $t\in[0,h_r),$ then $y(t)=e^{At}y(0)$ for $t\in[0,h_r).$ By using $u(t)=0$ for $t\in[-h_r,0],$ it is
obtained that $y(0)=x(0)$ from (\ref{i3}). Combining with the proof of Lemma \ref{lem1} and (\ref{x2}), the result follows.
So we omit the details.

Next, we consider the optimization problem with respect to the admissible control set set (\ref{i7}).

\begin{lemma}\label{lem4}
If a given linear feedback control $u(t) =K(t)\hat{y}(t|t)$ is the unique optimal solution for the problem of minimizing $J_T$ s.t (\ref{i4}), then $K(t)$ obeys the equations (\ref{6})-(\ref{e15}) with $\Omega(t)>0$.
\end{lemma}

Proof. The proof is presented in Appendix \ref{ap4}.

We now give the necessary and sufficient condition for the existence and uniqueness of the solution to the finite-horizon optimization problem.
\begin{theorem}\label{the1}
The problem of minimizing (\ref{2}) subject to (\ref{i4}) within the admissible control set (\ref{i7})
has a unique solution if and only if (\ref{x1})-(\ref{x5}) admits a solution such that the matrix
$\Omega(t)$ is strictly positive definite. The optimal control is as (\ref{contr}) and the optimal cost is given by (\ref{zz18}).
\end{theorem}

Proof. Combining with Lemmas \ref{lem2}-\ref{lem4}, the result follows directly.

\subsection{Finite-horizon optimal control problem of pseudo delay-free system with discounted cost function}

In this subsection ,we study the finite-horizon optimization problem of minimizing the discounted cost function subject to (\ref{i4}):
\begin{eqnarray}
J_T^{\alpha}&=&E\Big[
\int_{0}^{T }e^{-\alpha t}\Big(y'(t)Qy(t)+u(t)'Ru(t)\Big)dt\Big].\label{d3}
\end{eqnarray}
The discounted
setting is popular in many areas, such as in dynamic programming, reinforcement learning, and planning algorithms
for optimal control. See \cite{lavalle}, \cite{sutton} and references therein.

To solve the discounted LQR problem, we define the modified Riccati equation:
\begin{eqnarray}
-\frac{d}{dt}\hat{P}_{\alpha}(t)&=&\hat{P}_{\alpha}(t)A+A'\hat{P}_{\alpha}(t)+\alpha \hat{P}_{\alpha}(t)+Q\nonumber\\
&&-\Pi_{\alpha}(t,t),\label{d4}\\
P_{\alpha}(t)
&=&\hat{P}_{\alpha}(t)+\int_0^{h_r}e^{(A+\frac{\alpha}{2}I)'\theta}\Pi_{\alpha}(t+\theta,t+\theta)\nonumber\\
&&\times e^{(A+\frac{\alpha}{2}I)\theta}d\theta
,\label{d5}
\end{eqnarray}
where
\begin{eqnarray}
\Pi_{\alpha}(t,t)&=&K_{\alpha}'(t) \Omega_{\alpha}(t) K_{\alpha}(t), \nonumber\\
\Omega_{\alpha}(t) &=&R+\sum_{i=0}^{r}\bar{B}_i'e^{-A'h_i}P_{\alpha}(t)e^{-Ah_i}\bar{B}_i,\nonumber\\
K_{\alpha}(t)&=&\Omega_{\alpha}^{-1}(t)B'\hat{P}_{\alpha}(t),\nonumber
\end{eqnarray}
with $\hat{P}_{\alpha}(T)=0$ and $P_{\alpha}(T)=0.$
Following similar discussions to Lemma \ref{lem2} and Remark 5 in \cite{hszhang1}, the following result is in force.
\begin{lemma}\label{lem7}
The equation (\ref{d4})-(\ref{d5}) is equivalent to the following equations:
\begin{eqnarray}
-\dot{P}_{\alpha}(t)&=&P_{\alpha}(t)A+A'P_{\alpha}(t)+\alpha P_{\alpha}(t)+Q\nonumber\\
&&-e^{A'h_r}\Pi_{\alpha}(t+h_r,t+h_r)e^{Ah_r},\label{d6}
\end{eqnarray}
while $\Pi_{\alpha}(t+h_r,t+h_r)$ is given by
\begin{eqnarray}
\Pi_{\alpha}(t,t)&=&K_{\alpha}'(t)\Omega_{\alpha}(t)K_{\alpha}(t),\nonumber\\
\Omega_{\alpha}(t)&=&R+\sum_{i=0}^{r}\bar{B}_i'e^{-A'h_i}P_{\alpha}(t)e^{-Ah_i}\bar{B}_i,\nonumber\\
K_{\alpha}(t)&=&-\Omega_{\alpha}^{-1}(t)\Big[B'P_{\alpha}(t)
-B'\int_{0}^{h_r}e^{A'\theta}\Pi_{\alpha}(t+\theta,t+\theta)\nonumber\\
&&\times e^{A\theta}d\theta\Big],\nonumber
\end{eqnarray}
with terminal values $P_{\alpha}(T)=0$ and $\Pi_{\alpha}(T,T+\theta)=0$ for $\theta\in (0,h_r].$
\end{lemma}

It is now in the position to give the solution to the discounted LQR problem.
\begin{theorem}\label{theorem}
The problem of minimizing (\ref{d3}) subject to (\ref{i4}) within the admissible control set (\ref{i7})
has a unique solution if and only if (\ref{d4})-(\ref{d5}) admits a solution such that the matrix
$\Omega_{\alpha}(t)$ is strictly positive definite. The optimal control is as
\begin{eqnarray}
u(t)=K_{\alpha}(t)\hat{y}(t|t),\label{d7}
\end{eqnarray} and the optimal cost is given by
\begin{eqnarray}
J_T^{\ast}&=&E\Big(y'(0)P_{\alpha}(0)y(0)-y'(0)\int_{0}^{h_r}
\Pi_{\alpha}(0,\theta)\hat{y}(0|\theta)d\theta\Big).\label{d8}
\end{eqnarray}
\end{theorem}

Proof. The proof is presented in Appendix \ref{apd}.

\section{Solution to the \emph{Problem}}
Based on the above results for the finite-horizon optimization problem, we discuss the mean-square stabilization problem.
Sufficient and necessary conditions are to be derived for the exponential mean-square stabilization of system (\ref{i1}). The key is to
investigate the properties of the modified Riccati equations (\ref{x1})-(\ref{x5}) and (\ref{d4})-(\ref{d5}) when the time $t$ tends to $-\infty.$
Firstly, we give the necessary condition for the mean-square stabilization for system (\ref{i1}).
\begin{theorem}\label{lem5}
Assume that the system (\ref{i1}) is exponentially mean-square stabilizable in the sense of Definition \ref{d2},
then the following modified algebraic Riccati equation (\ref{x6})-(\ref{x10}) has a solution ${P}\geq \hat{P}>0$,

\begin{eqnarray}
  0&=&A'\hat{P}+\hat{P}A-\Pi(0)+I,\label{x6}\\
  {P}&=&\hat{P}+\int_0^{h_r} e^{A'\theta}\Pi(0) e^{A\theta}d\theta,\label{x7}
  \end{eqnarray}
where
 \begin{eqnarray}
   \Pi(0)&=&K' \Omega K, \label{x8} \\
  \Omega &=&I+\sum_{i=0}^{r}\bar{B}_i'e^{-A'h_i}P(t)e^{-Ah_i}\bar{B}_i,\label{x9} \\
K&=&\Omega^{-1}B'\hat{P}.\label{x10}
  \end{eqnarray}

\end{theorem}

Proof. The proof is put in Appendix \ref{ap5}.

We then present the sufficient condition for the exponential mean-square stabilization by defining a new Lyapunov function.
\begin{theorem}\label{lem6}
Assume that the following equation has a unique solution ${P}_{\alpha}\geq \hat{P}_{\alpha}>0,$
\begin{eqnarray}
  0&=&A'\hat{P}_{\alpha}+\hat{P}_{\alpha}A+\alpha P_{\alpha}-\Pi_{\alpha}(0)+I,\label{ex6}\\
  {P}_{\alpha}&=&\hat{P}_{\alpha}+\int_0^{h_r} e^{(A+\frac{\alpha }{2}I)'\theta}\Pi_{\alpha}(0) e^{(A+\frac{\alpha }{2}I)\theta}d\theta,\label{ex7}
  \end{eqnarray}
where
 \begin{eqnarray}
   \Pi_{\alpha}(0)&=&K_{\alpha}' \Omega_{\alpha} K_{\alpha}, \label{ex8} \\
  \Omega_{\alpha} &=&I+\sum_{i=0}^{r}\bar{B}_i'e^{-A'h_i}P_{\alpha}(t)e^{-Ah_i}\bar{B}_i,\label{ex9} \\
K_{\alpha}&=&\Omega_{\alpha}^{-1}B'\hat{P}_{\alpha},\label{ex10}
  \end{eqnarray}
then the system (\ref{i1}) is exponentially mean-square stable with the controller $u(t)=K_{\alpha}\hat{y}(t|t)$ where $K_{\alpha}$ is given by (\ref{ex10}).
\end{theorem}

Proof. The proof is formulated in Appendix \ref{ap6}.

\section{Conclusions}

This paper studied the stabilization problem for the It\^{o} systems with both multiplicative noise and multiple delays.
Sufficient and necessary conditions have been obtained for the exponential mean-square stabilization in terms of modified
Riccati equations. The main technique is to reduce the original system with multiple delays to the pseudo delay-free one and study the
finite-horizon optimization problems for the pseudo system with standard and discounted linear quadratic cost functions.

\appendix

\section{Proof of Lemma \ref{lem1}}\label{ap2}
Using Lemma \ref{lem2}, the equations (\ref{6})-(\ref{e15}) admit a solution such that the matrix
$\Omega(t)>0$. Applying It\^{o}'s formula to $y'(t)\big[P(t)y(t)-\int_{0}^{h_r}
\Pi(t,t+\theta)\hat{y}(t|t+\theta)d\theta\big]$ and combining with the equations (\ref{6})-(\ref{e15}),
we have
\begin{eqnarray}
&&
d\Big\{y'(t)\Big[P(t)y(t)-\int_{0}^{h_r}\Pi(t,t+\theta)
\hat{y}(t|t+\theta)d\theta\Big]\Big\}\nonumber\\
&=&\Big\{\Big(Ay(t)+\sum_{i=0}^{r}e^{-Ah_i}B_iu(t)\Big)'\Big[P(t)y(t)\nonumber\\
&&-\int_{0}^{h_r}\Pi(t,t+\theta)
\hat{y}(t|t+\theta)d\theta\Big]+y'(t)\dot{P}(t)y(t)\nonumber\\
&&+y'(t)P(t)\Big(Ay(t)+\sum_{i=0}^{r}e^{-Ah_i}B_iu(t)\Big)\nonumber\\
&&+u'(t)\sum_{i=0}^{r}\bar{B}_i'e^{-A'h_i}P(t)e^{-Ah_i}\bar{B}_iu(t)\nonumber\\
&&-y'(t)\Pi(t,t+h_r)y(t)+y'(t)\Pi(t,t)\hat{y}(t|t)\nonumber\\
&&-y'(t)\int_{0}^{h_r}\frac{\partial}{\partial t}\Pi(t,t+\theta)
\hat{y}(t|t+\theta)d\theta-y'(t)\nonumber\\
&&\times \int_{0}^{h_r}\Pi(t,t+\theta)
\Big(A\hat{y}(t|t+\theta)+\sum_{i=0}^{r}e^{-Ah_i}B_iu(t)\Big)d\theta\Big\}dt\nonumber\\
&&+\Big\{\Big[\sum_{i=0}^{r}e^{-Ah_i}\bar{B}_iu(t)\Big]'P(t)y(t)+y'(t)P(t)\nonumber\\
&&\times \Big[\sum_{i=0}^{r}e^{-Ah_i}\bar{B}_iu(t)\Big]-\Big[\sum_{i=0}^{r}e^{-Ah_i}\bar{B}_iu(t)\Big]'\nonumber\\
&&\times \int_{t}^{t+h_r}
\Pi(t,\theta)\hat{y}(t|\theta)d\theta \Big\}dw(t)\nonumber\\
&=&\Big[-y'(t)Qy(t)+y'(t)\Pi(t,t)
\hat{y}(t|t)+2u'(t)B'P(t)y(t)\nonumber\\
&&
+u'(t)\sum_{i=0}^{r}\bar{B}_i'e^{-A'h_i}P(t)e^{-Ah_i}\bar{B}_iu(t)\nonumber\\
&&-u'(t)B'\int_t^{t+h_r}\Pi(t,\theta)\hat{y}(t|\theta)d\theta\nonumber\\
&&-y'(t)\int_t^{t+h_r}\Pi(t,\theta)d\theta Bu(t)\Big]dt\nonumber\\
&&+\Big\{\Big[\sum_{i=0}^{r}e^{-Ah_i}\bar{B}_iu(t)\Big]'P(t)y(t)\nonumber\\
&&+y'(t)P(t)\Big[\sum_{i=0}^{r}e^{-Ah_i}\bar{B}_iu(t)\Big]\nonumber\\
&&-\Big[\sum_{i=0}^{r}e^{-Ah_i}\bar{B}_iu(t)\Big]'\int_{t}^{t+h_r}
\Pi(t,\theta)\hat{y}(t|\theta)d\theta \Big\}dw(t).\label{zz13}
\end{eqnarray}
Taking integral from $0$ to $T$ on both sides of (\ref{zz13}) and then taking expectation, we have
\begin{eqnarray}
J_T
&=&E\Big(y'(0)P(0)y(0)-y'(0)\int_{0}^{h_r}
\Pi(0,\theta)\hat{y}(0|\theta)d\theta\Big)\nonumber\\
&&+E\int_{0}^T\Big(u'(t)\Omega(t)u(t)-2u'(t)\Omega(t) K(t) y(t)\nonumber\\
&&+y'(t)\Pi(t,t)\hat{y}(t|t)\Big)dt\nonumber\\
&=&E\Big(y'(0)P(0)y(0)-y'(0)\int_{0}^{h_r}
\Pi(0,\theta)\hat{y}(0|\theta)d\theta\Big)\nonumber\\
&&+E\int_{0}^T\Big(u(t)-K(t)\hat{y}(t|t)\Big)'\Omega(t)\Big(u(t)\nonumber\\
&&-K(t)\hat{y}(t|t)\Big)dt,\label{13}
\end{eqnarray}
where the fact of $E\Big\{\Big[y(t)-\hat{y}(t|t)\Big]'\hat{y}(t|t)\Big\}=0$ has been used in the derivation of the above equality. Note that $\Omega(t)>0$, the optimal control exists uniquely. Furthermore, the optimal control (\ref{contr}) and cost function (\ref{zz18}) follows from (\ref{13}) directly combining with Lemma \ref{lem2}.

\section{Proof of Lemma \ref{lem4}}\label{ap4}

Consider the optimization problem for the controller set $\{u(t):u(t)=K(t)\hat{y}(t|t)\}$ with respect to the matrix $K(t)$. The cost function is
\begin{eqnarray}
J_T&=&E\Big[\int_0^T\big[y'(t)Qy(t)+\hat{y}'(t|t)K'(t)RK(t)\hat{y}(t|t)\big]dt\nonumber\\
&&+y'(T)Hy(T)\Big]\nonumber\\
&=&tr\Big[\int_0^T\big[QY(t)+K'(t)RK(t)\hat{Y}(t|t)\big]dt+HY(T)\Big],\nonumber\\
\end{eqnarray}
where $Y(t)=E[y(t)y'(t)]$ and $\hat{Y}(t|t)=E[\hat{y}(t|t)\hat{y}'(t|t)].$
The system under the controller $u(t)=K(t)\hat{y}(t|t)$ is reduced to
\begin{eqnarray}
dy(t)&=&\big[Ay(t)+BK(t)\hat{y}(t|t)\big]dt+\nonumber\\
&&\sum_{i=0}^{r}e^{-Ah_i}\bar{B}_iK(t)\hat{y}(t|t)dw_i(t+h_i).
\end{eqnarray}
In this case,
\begin{eqnarray}
&&d[y(t)y'(t)]\nonumber\\
&=&\big[Ay(t)+BK(t)\hat{y}(t|t)\big]y'(t)dt\nonumber\\
&&+\sum_{i=0}^{r}e^{-Ah_i}\bar{B}_iK(t)\hat{y}(t|t)y'(t)dw_i(t+h_i)\nonumber\\
&&+y(t)\big[Ay(t)+BK(t)\hat{y}(t|t)\big]'dt\nonumber\\
&&+y(t)\sum_{i=0}^{r}\Big(e^{-Ah_i}\bar{B}_iK(t)\hat{y}(t|t)\Big)'dw_i(t+h_i)\nonumber\\
&&
+\sum_{i=0}^{r}e^{-Ah_i}\bar{B}_iK(t)\hat{y}(t|t)\hat{y}'(t|t)K'(t)\bar{B}_i'e^{-A'h_i}dt,\nonumber
\end{eqnarray}
that is,
\begin{eqnarray}
\frac{d}{dt}Y(t)
&=&AY(t)+BK(t)\hat{Y}(t|t)+Y(t)A'+\hat{Y}(t|t)\nonumber\\
&&\times K'(t)B'+\sum_{i=0}^{r}e^{-Ah_i}\bar{B}_iK(t)\hat{Y}(t|t)K'(t)\nonumber\\
&&\times \bar{B}_i'e^{-A'h_i}.
\end{eqnarray}
In addition, it is obtained that
\begin{eqnarray}
\frac{\partial}{\partial t}\hat{Y}(t|\theta)
&=&A\hat{Y}(t|\theta)+BK(t)\hat{Y}(t|t)+\hat{Y}(t|\theta)A'\nonumber\\
&&+\hat{Y}(t|t)K'(t)B',\nonumber
\end{eqnarray}
thus, we have
\begin{eqnarray}
&&\frac{d}{dt}\int_t^{t+h_r}\hat{Y}(t|\theta)\Pi'(t,\theta)d\theta\nonumber\\
&=&Y(t)\Pi'(t,t+h_r)-\hat{Y}(t|t)\Pi'(t,t)+\int_t^{t+h_r}\frac{\partial}{\partial t}\hat{Y}(t|\theta)\nonumber\\
&&\times \Pi'(t,\theta)d\theta+\int_t^{t+h_r}\hat{Y}(t|\theta)\frac{\partial}{\partial t}\Pi'(t,\theta)d\theta.
\end{eqnarray}
Using the Lagrange multiplier approach, the cost function can be reformulated as follows:
\begin{eqnarray}
J_T&=&\int_0^Ttr\Big[QY(t)+K'(t)RK(t)\hat{Y}(t|t)+\big[AY(t)\nonumber\\
&&+BK(t)\hat{Y}(t|t)+Y(t)A'+\hat{Y}(t|t) K'(t)B'\nonumber\\
&&+\sum_{i=0}^{r}e^{-Ah_i}\bar{B}_iK(t)\hat{Y}(t|t)K'(t)\bar{B}_i'e^{-A'h_i}-\dot{Y}(t)\big]\nonumber\\
&&\times P'(t)-\int_t^{t+h_r}\big[A\hat{Y}(t|\theta)+BK(t)\hat{Y}(t|t)
+\hat{Y}(t|\theta)A'\nonumber\\
&&+\hat{Y}(t|t)K'(t)B'-\frac{\partial}{\partial t}\hat{Y}(t|\theta)\big]\Pi'(t,\theta)d\theta\Big]dt\nonumber\\
&&+tr[HY(T)],\nonumber
\end{eqnarray}
where $P(\cdot),\Pi(\cdot,\cdot)$ are matrix parameters with compatible dimension.
By making some algebraic transformation, it is further rewritten as
\begin{eqnarray}
J_T&=&\int_0^Ttr\Big[QY(t)+K'(t)RK(t)\hat{Y}(t|t)+\big[AY(t)\nonumber\\
&&+BK(t)\hat{Y}(t|t)+Y(t)A'+\hat{Y}(t|t)K'(t)B'\nonumber\\
&&
+\sum_{i=0}^{r}e^{-Ah_i}\bar{B}_iK(t)\hat{Y}(t|t)K'(t)\bar{B}_i'e^{-A'h_i}\big]P'(t)\nonumber\\
&&+Y(t)\dot{P}'(t)-\int_t^{t+h_r}\big[A\hat{Y}(t|\theta)+BK(t)\hat{Y}(t|t)\nonumber\\
&&+\hat{Y}(t|\theta)A'+\hat{Y}(t|t)K'(t)B'-\frac{\partial}{\partial t}\hat{Y}(t|\theta)\big]\nonumber\\
&&\times \Pi'(t,\theta)d\theta\Big]dt+tr[HY(T)]-P(T)Y(T)\nonumber\\
&&+P(0)Y(0)\nonumber\\
&=&\int_0^Ttr\Big[QY(t)+K'(t)RK(t)\hat{Y}(t|t)+\big[AY(t)\nonumber\\
&&+BK(t)\hat{Y}(t|t)+Y(t)A'+\hat{Y}(t|t)K'(t)B'\nonumber\\
&&+\sum_{i=0}^{r}e^{-Ah_i}\bar{B}_iK(t)\hat{Y}(t|t)K'(t)\bar{B}_i'e^{-A'h_i}\big]P'(t)\nonumber\\
&&+Y(t)\dot{P}'(t)-\int_t^{t+h_r}\big[A\hat{Y}(t|\theta)+BK(t)\hat{Y}(t|t)\nonumber\\
&&
+\hat{Y}(t|\theta)A'+\hat{Y}(t|t)K'(t)B'\big]\Pi'(t,\theta)d\theta\nonumber\\
&&-Y(t)\Pi'(t,t+h_r)+\hat{Y}(t|t)\Pi'(t,t)\nonumber\\
&&-\int_t^{t+h_r}\hat{Y}(t|\theta)\frac{\partial}{\partial t}\Pi'(t,\theta)d\theta\Big]dt+trHY(T)\nonumber\\
&&-P(T)Y(T)+P(0)Y(0)+\int_T^{T+h_r}\hat{Y}(T|\theta)\Pi'(T,\theta)d\theta\nonumber\\
&&-\int_0^{h_r}\hat{Y}(0|\theta)\Pi'(0,\theta)d\theta.
\end{eqnarray}
Taking partial differential yields that
\begin{eqnarray}
0&=&\frac{\partial J_T}{\partial Y(t)}=Q+A'P(t)+P(t)A+\dot{P}(t)\nonumber\\
&&~~~~~~~~~~~~~-\Pi(t,t+h_r),\nonumber\\
0&=&\frac{\partial J_T}{\partial \hat{Y}(t|\theta-h)}=-\frac{\partial}{\partial t} \Pi(t,\theta)-A'\Pi(t,\theta)-\Pi(t,\theta)A,\nonumber\\
0&=&\frac{\partial J_T}{\partial \hat{Y}(t|t)}\nonumber\\
&=&K'(t)RK(t)+K'(t)B'P(t)+P(t)BK(t)\nonumber\\
&&+\sum_{i=0}^{r}K'(t)\bar{B}_i'e^{-A'h_i}P(t)e^{-Ah_i}\bar{B}_iK(t)\nonumber\\
&&-K'(t)B'\int_t^{t+h_r}\Pi(t,\theta)d\theta\nonumber\\
&&-\int_t^{t+h_r}\Pi(t,\theta)d\theta BK(t)+\Pi(t,t)\nonumber\\
&=&K'(t)\Omega(t) K(t)+K'(t)\Big(B'P(t)-B'\nonumber\\
&&\times \int_t^{t+h_r}\Pi(t,\theta)d\theta\Big)+\Big(P(t)B-\int_t^{t+h_r}\Pi(t,\theta)d\theta B\Big)\nonumber\\
&&\times K(t)+\Pi(t,t),\nonumber\\
0&=&\frac{\partial J_T}{\partial K(t)}\nonumber\\
&=&RK(t)\hat{Y}(t|t)+RK(t)\hat{Y}'(t|t)+B'P(t)\hat{Y}'(t|t)\nonumber\\
&&+B'P'(t)\hat{Y}(t|t)+\sum_{i=0}^{r}\bar{B}_i'e^{-A'h_i}P(t)e^{-Ah_i}\bar{B}_iK(t)\nonumber\\
&&\times \hat{Y}'(t|t)-\int_t^{t+h_r}B'\Pi(t,\theta)d\theta \hat{Y}'(t|t)\nonumber\\
&&+\sum_{i=0}^{r}\bar{B}_i'e^{-A'h_i}P'(t)e^{-Ah_i}\bar{B}_iK(t)\hat{Y}(t|t)\nonumber\\
&&-\int_t^{t+h_r}B'\Pi'(t,\theta)d\theta \hat{Y}(t|t)\nonumber\\
&=&[\Omega(t) K(t)+B'P'(t)-\int_t^{t+h_r}B'\Pi'(t,\theta)d\theta]\hat{Y}(t|t)\nonumber\\
&&+[\Omega(t) K(t)+B'P(t)-\int_t^{t+h_r}B'\Pi(t,\theta)d\theta]\hat{Y}'(t|t),\nonumber
\end{eqnarray}
with $P(T)=H$ and $\Pi(T,\theta)=0$.
Thus, we have the following equation:
\begin{eqnarray}
-\dot{P}(t)&=&Q+A'P(t)+P(t)A-\Pi(t,t+h_r),\nonumber\\
-\frac{\partial}{\partial t} \Pi(t,\theta)&=&A'\Pi(t,\theta)+\Pi(t,\theta)A,~~~\Pi(T,\theta)=0,\nonumber\\
\Pi(t,t)&=&K'(t)\Omega(t) K(t),\nonumber\\
0&=&\Omega(t) K(t)+B'P(t)-\int_t^{t+h_r}B'\Pi(t,\theta)d\theta.\nonumber
\end{eqnarray}
Using the unique existence of the optimal controller, we have the positive definiteness of the matrix $\Omega(t)>0$. Thus, (\ref{6})-(\ref{e15}) admits a solution with $\Omega(t)>0$.

\section{Proof of Theorem \ref{theorem}}\label{apd}

``Necessity" By applying similar procedures to Lemma \ref{lem4}, the necessity follows directly. To avoid duplication, we omit the details.

``Sufficiency" Using Lemma \ref{lem7}, the equation (\ref{d6}) admits a solution such that the matrix
$\Omega_{\alpha}(t)>0$. Applying It\^{o}'s formula to $e^{\alpha t}y'(t)\big[P_{\alpha}(t)y(t)-\int_{0}^{h_r}
\Pi_{\alpha}(t,t+\theta)\hat{y}(t|t+\theta)d\theta\big]$ and combining with the equations (\ref{d4})-(\ref{d5}),
we have
\begin{eqnarray}
&&d\Big[e^{\alpha t}y'(t)\Big(P_{\alpha}(t)y(t)-\int_{0}^{h_r}\Pi_{\alpha}(t,t+\theta)
\hat{y}(t|t+\theta)d\theta\Big)\Big]\nonumber\\
&=&e^{\alpha t}\Big\{\alpha y'(t)\Big(P_{\alpha}(t)y(t)-\int_{0}^{h_r}\Pi_{\alpha}(t,t+\theta)
\hat{y}(t|t+\theta)d\theta\Big)\nonumber\\
&&+ \Big(Ay(t)+\sum_{i=0}^{r}e^{-Ah_i}B_iu(t)\Big)'\Big[P(t)y(t)\nonumber\\
&&-\int_{0}^{h_r}\Pi(t,t+\theta)
\hat{y}(t|t+\theta)d\theta\Big]+y'(t)\dot{P}(t)y(t)\nonumber\\
&&+y'(t)P(t)\Big(Ay(t)+\sum_{i=0}^{r}e^{-Ah_i}B_iu(t)\Big)\nonumber\\
&&+u'(t)\sum_{i=0}^{r}\bar{B}_i'e^{-A'h_i}P(t)e^{-Ah_i}\bar{B}_iu(t)\nonumber\\
&&-y'(t)\Pi(t,t+h_r)y(t)+y'(t)\Pi(t,t)\hat{y}(t|t)\nonumber\\
&&-y'(t)\int_{0}^{h_r}\frac{\partial}{\partial t}\Pi(t,t+\theta)
\hat{y}(t|t+\theta)d\theta\nonumber\\
&&-y'(t)\int_{0}^{h_r}\Pi(t,t+\theta)
\Big(A\hat{y}(t|t+\theta)+\sum_{i=0}^{r}e^{-Ah_i}\nonumber\\
&&\times B_iu(t)\Big)d\theta\Big\}dt+e^{\alpha t}\Big\{\Big[\sum_{i=0}^{r}e^{-Ah_i}\bar{B}_iu(t)\Big]'P(t)y(t)\nonumber\\
&&+y'(t)P(t)\Big[\sum_{i=0}^{r}e^{-Ah_i}\bar{B}_i u(t)\Big]\nonumber\\
&&-\Big[\sum_{i=0}^{r}e^{-Ah_i}\bar{B}_iu(t)\Big]'\int_{t}^{t+h_r}
\Pi(t,\theta)\hat{y}(t|\theta)d\theta \Big\}dw(t)\nonumber\\
&=&e^{\alpha t}\Big[-y'(t)Qy(t)+y'(t)\Pi_{\alpha}(t,t)
\hat{y}(t|t)\nonumber\\
&&+2u'(t)B'P_{\alpha}(t)y(t)
+u'(t)\sum_{i=0}^{r}\bar{B}_i'e^{-A'h_i}P_{\alpha}(t)e^{-Ah_i}\nonumber\\
&&\times \bar{B}_iu(t)-u'(t)B'\int_t^{t+h_r}\Pi_{\alpha}(t,\theta)\hat{y}(t|\theta)d\theta\nonumber\\
&&
-y'(t) \int_t^{t+h_r}\Pi_{\alpha}(t,\theta)d\theta Bu(t)\Big]dt+\Big\{\Big[\sum_{i=0}^{r}e^{-Ah_i}\bar{B}_iu(t)\Big]'\nonumber\\
&&\times P_{\alpha}(t)y(t)
+y'(t)P_{\alpha}(t)\Big[\sum_{i=0}^{r}e^{-Ah_i}\bar{B}_iu(t)\Big]\nonumber\\
&&-\Big[\sum_{i=0}^{r}e^{-Ah_i}\bar{B}_iu(t)\Big]'\int_{t}^{t+h_r}
\Pi_{\alpha}(t,\theta)\hat{y}(t|\theta)d\theta \Big\}dw(t).\nonumber
\end{eqnarray}
Taking integral from $0$ to $T$ and then taking expectation on both sides of the above equation, we have
\begin{eqnarray}
J_T^{\alpha}
&=&E\Big(y'(0)P_{\alpha}(0)y(0)-y'(0)\int_{0}^{h_r}
\Pi_{\alpha}(0,\theta)\hat{y}(0|\theta)d\theta\Big)\nonumber\\
&&+E\int_{0}^Te^{\alpha t}\Big(u'(t)\Omega_{\alpha}(t)u(t)-2u'(t)\Omega_{\alpha}(t) K_{\alpha}(t) y(t)\nonumber\\
&&+y'(t)\Pi_{\alpha}(t,t)\hat{y}(t|t)\Big)dt\nonumber\\
&=&E\Big(y'(0)P_{\alpha}(0)y(0)-y'(0)\int_{0}^{h_r}
\Pi_{\alpha}(0,\theta)\hat{y}(0|\theta)d\theta\Big)\nonumber\\
&&+E\int_{0}^Te^{\alpha t}\Big(u(t)-K_{\alpha}(t)\hat{y}(t|t)\Big)'\Omega_{\alpha}(t)\Big(u(t)\nonumber\\
&&-K_{\alpha}(t)\hat{y}(t|t)\Big)dt,\label{d9}
\end{eqnarray}
where the fact of $E\Big\{\Big[y(t)-\hat{y}(t|t)\Big]'\hat{y}(t|t)\Big\}=0$ has been used in the derivation of the above equality. Note that $\Omega_{\alpha}(t)>0$, the optimal control exists uniquely. Furthermore, the optimal control (\ref{d7}) and optimal cost function (\ref{d8}) follows from (\ref{d9}) directly.

\section{Proof of Theorem \ref{lem5}}\label{ap5}

In view of Theorem \ref{the1}, the fact that $R=I>0$ can ensure the existence of the solution to (\ref{x1})-(\ref{x5}) with $\Omega(t)>0$. Re-denote the solution $P(t),\hat{P}(t)$ and
$\Pi(t,t+\theta)$ of (\ref{x1})-(\ref{x5}) as $P_T(t),\hat{P}_T(t)$ and $\Pi_T(t,t+\theta)$ respectively, with
the terminal time $T$ and the terminal values $P(T)=H=0,\hat{P}(T)=0$ and
$\Pi(T,T+\theta)=0.$ We first show that $\hat{P}_T(t)$ of (\ref{x1}) and $P_T(t)$  of (\ref{x2}) are convergent.
Based on Corollary \ref{cor2}, the optimal cost becomes
${J_T}^*=E\Big(x_0'\hat{P}_T(0)x_0\Big).$
Noting the time-invariance of (\ref{6})-(\ref{e15}) with respect to $T$, i.e., for $ t\leq T,$
\begin{eqnarray}
P_T(t)=P_{T-t}(0),
\Pi_T(t,t+\theta)=\Pi_{T-t}(0,\theta),
\theta\in[0,h_r].\nonumber
\end{eqnarray}
Thus, for any $T_1>T>t$ and for all $ x_0\neq 0,$ we have
\begin{eqnarray}
&&x_0'\hat{P}_{T_1}(t)x_0=x_0'\hat{P}_{T_1-t}(0)x_0={J_{T_1-t}}^* \nonumber\\
&\geq&
{J_{T-t}}^*=x_0'\hat{P}_{T-t}(0)x_0=x_0'\hat{P}_{T}(t)x_0.\nonumber
\end{eqnarray}
Since $x_0$ is arbitrary, thus $\hat{P}_{T_1}(t)\geq \hat{P}_{T}(t).$ Similarly, if $ t_1< t_2\leq T,$
\begin{eqnarray}
&&x_0'\hat{P}_{T}(t_1)x_0=x_0'\hat{P}_{T-t_1}(0)x_0={J_{T-t_1}}^* \nonumber\\
&\geq&
{J_{T-t_2}}^*=x_0'\hat{P}_{T-t_2}(0)x_0=x_0'\hat{P}_{T}(t_2)x_0.\nonumber
\end{eqnarray}
That is, $\hat{P}_{T}(t_1)\geq \hat{P}_{T}(t_2).$ Thus, $\hat{P}_T(t)$ is monotonically increasing with respect to $T$
and is monotonically decreasing with respect to $t.$

We then show the uniform boundedness of $\hat{P}_T(t)$. Since system (\ref{i1}) is exponentially stabilizable in the sense of
Definition \ref{d2}, together with (\ref{i3}), there exists a positive constant $\delta$ such that
\begin{eqnarray}
&&e^{\alpha t}E\|y(t)\|^2\nonumber\\
&\leq& \delta e^{\alpha t}\Big(E\|x(t)\|^2+\sum_{i=1}^{r}\int_t^{t+h_i}\|e^{A(t-s)}B_iu(s-h_i)\|^2ds\nonumber\\
&&+\sum_{i=1}^{r}\int_t^{t+h_i}\|e^{A(t-s)}\bar{B}_iu(s-h_i)\|^2ds\Big)\nonumber\\
&\rightarrow&0,~~~ t\rightarrow \infty,\nonumber
\end{eqnarray}
where the last limit holds for $\lim_{t\rightarrow\infty}e^{\alpha t}E\|x(t)\|^2=0$ and
$\lim_{t\rightarrow\infty}e^{\alpha t}E\|u(t)\|^2=0.$
Together with the exponential stability of $u(t),$ we have the boundness of the cost function $J_T$ under the stabilizing controller. In fact, there exists
a positive constant $\mu$ such that $e^{\alpha t}E\|y(t)\|^2\leq \mu\|x_0\|^2$ and
$e^{\alpha t}E\|u(t)\|^2\leq \mu\|x_0\|^2.$ This further implies that there exists a positive constant $\beta$ such that
\begin{eqnarray}
E\int_0^\infty\Big(y'(t)Qy(t)+u'(t)Ru(t)\Big)dt\leq \beta\|x_0\|^2.\nonumber
\end{eqnarray}
Thus
\begin{eqnarray}
J_T^*=x_0'\hat{P}_T(0)x_0<\beta\|x_0\|^2,\nonumber
\end{eqnarray}
that is, $\hat{P}_T(0)$ is uniformly bounded.
Recalling the monotonicity of $\hat{P}_T(t)$, it yields that $\hat{P}_T(t)$ is convergent, i.e.,
\begin{eqnarray}
\lim_{t\rightarrow -\infty} \hat{P}_T(t)=\lim_{t\rightarrow -\infty}
\hat{P}_{T-t}(0)=\lim_{T\rightarrow \infty} \hat{P}_{T}(0)\doteq \hat{P},\nonumber
\end{eqnarray}
where $\hat{P}$ is a constant matrix which is independent of
$t$.

Consider the optimal cost (\ref{i9}) in Corollary \ref{cor1}, we have that $P_T(t)$ is
monotonically increasing with respect to $T$ and is monotonically decreasing with respect to $t.$ Moreover,
$P_T(0)$ is uniformly bounded. The discussion is similar to that of $\hat{P}_T(t)$, so we omit the details.
This implies that $P_T(t)$ is convergent, i.e.,
\begin{eqnarray}\lim_{t\rightarrow -\infty} P_T(t)=\lim_{t\rightarrow -\infty}
P_{T-t}(0)=\lim_{T\rightarrow \infty} P_{T}(0)\doteq P,\nonumber
\end{eqnarray}
where $P$ is a constant matrix which is independent of
$t$.
Let $t\rightarrow-\infty$ in the equations (\ref{x1})-(\ref{x5}), we immediately have (\ref{x6})-(\ref{x10}).

Secondly, we show the strictly positive definiteness of
the matrix $\hat{P}$.  Otherwise, there exists $z\neq0,$ such that
$z'\hat{P}z=0.$
Similar to (\ref{zz13}) and (\ref{13}), by applying It\^{o}'s formula to $y'(t)[Py(t)-\int_{0}^{h_r}\Pi(\theta)\hat{y}(t|t+\theta)d\theta]$
where $\Pi(\theta)=e^{A'\theta}\Pi(0)e^{A\theta}$, $P$ and $\Pi(\theta)$ are as in (\ref{x7})-(\ref{x10}),  it follows that
\begin{eqnarray}
&&E\Big\{y'(T)\big[Py(T)-\int_{0}^{h_r}\Pi(\theta)\hat{y}(T|T+\theta)d\theta\big]\Big\}-Ex_0'\hat{P}x_0\nonumber\\
&=&E\int_0^T \big[-y'(t)y(t)
+y'(t)\Pi(0)\hat{y}(t|t)-2u'(t)\Omega K\hat{y}(t|t)\nonumber\\
&&+u'(t)\sum_{i=0}^{r}\bar{B}_i'e^{-A'h_i}P(t)e^{-Ah_i}\bar{B}_iu(t)\big]dt, \nonumber
\end{eqnarray}
Let $u(t)=K\hat{y}(t|t),\ t\in[0,T]$,  thus
\begin{eqnarray}
0&\leq &E\int_0^T\big[y'(t)y(t)+u(t)'u(t)\big]dt\nonumber \\
&=&
-E\Big[y'(T)\Big(Py(T)-\int_{0}^{h_r}\Pi(\theta)
\hat{y}(T|T+\theta)d\theta\Big)\Big]\nonumber\\
&&+x_0'\hat{P}x_0+E\int_0^T \Big(y'(t)\Pi(0)\hat{y}(t|t)-2u'(t)\Omega K y(t)\nonumber\\
&&+u'(t)\Omega u(t)\Big)dt\nonumber \\
&=&-E\Big[y'(T)\Big(Py(T)-\int_{0}^{h_r}\Pi(\theta)
\hat{y}(T|T+\theta)d\theta\Big)\Big]\nonumber\\
&&+x_0'\hat{P}x_0.\nonumber
\end{eqnarray}
Now let $x(0)=z$  where $z$ is given as $z'\hat{P}z=0$. Then, $x_0'\hat{P}x_0=0$. Thus
\begin{eqnarray}
0&\leq & E\Big(\int_0^Ty'(t)y(t)dt+\int_0^Tu(t)'u(t)dt\Big)\nonumber\\
&=&-E\Big(y'(T)Py(T)-y'(T)\int_0^{h_r}\Pi(\theta)\hat{y}(T|T+\theta)d\theta\Big).\nonumber\\\label{sz1}
\end{eqnarray}
Further note that $\Pi(\theta)\geq 0$ and $\hat{P}\geq 0$ as shown in the above, we have
\begin{eqnarray}
&&E\big[y'(T)Py(T)-y'(T)\int_0^{h_r}\Pi(\theta)\hat{y}(T|T+\theta)d\theta\big]\nonumber\\
&\geq &E\big[y'(T)Py(T)-y'(T)\int_0^{h_r}\Pi(\theta)\hat{y}(T|T+\theta)d\theta\nonumber\\
&&
-\int_0^{h_r}\tilde{y}'(T|T+\theta)
\Pi(\theta)\tilde{y}(T|T+\theta)d\theta\big]\nonumber\\
&=&E\big[y'(T)\hat{P}y(T)\big]\geq 0,\label{sz2}
\end{eqnarray}
where $\tilde{y}(t|t+\theta)=y(t)-\hat{y}(t|t+\theta)$ and $E[\tilde{y}'(T|T+\theta)\hat{y}(T|T+\theta)]=0$ have been used in the above.
Thus, it follows from (\ref{sz1}) and (\ref{sz2}) that
\begin{eqnarray}0\leq  E\int_0^T\big[y'(t)y(t)+u'(t)u(t)\big]dt
\leq 0.\nonumber
\end{eqnarray}
This implies that
\begin{eqnarray}
E[y'(t)y(t)]=0,
E[u'(t)u(t)]=0, t\geq 0.\nonumber
\end{eqnarray}
Then, it is obtained that $y(t)=0$ and $u(t)=0,\ t\geq 0$, $a.s.$. System (\ref{i4}) is thus
now reduced to
\begin{eqnarray}
dy(t)=Ay(t)dt,y(0)=x_0=z\neq0,\nonumber
\end{eqnarray}
with the output $y(t)=0$ a.s., this is a contradiction with
the observability of the system $(A,I)$. Thus, the matrix
$\hat{P}$ is positive definite. Together with (\ref{x7})-(\ref{x9}), $P\geq \hat{P}>0$ follows.
The proof is now completed.

\section{Proof of Theorem \ref{lem6}}\label{ap6}

We will prove that the system (\ref{i3}) is exponentially mean-square stabilizable under the controller $u(t)=K\hat{y}(t|t).$
Define the Lyapunov function candidate as
\begin{eqnarray}
V(t,y(t))&=&e^{\alpha t}E\big[y'(t)P_{\alpha}y(t)
-y'(t)\int_0^{h_r}\Pi_{\alpha}(\theta)\nonumber\\
&&\times \hat{y}(t|t+\theta)d\theta\big],t\geq 0,\label{kkv}
\end{eqnarray}
where $\Pi_{\alpha}(\theta)=e^{A'\theta}\Pi_{\alpha}(0)e^{A\theta}$. It is obvious that
\begin{eqnarray}
&&V(t,y(t))\nonumber\\
&\geq&e^{\alpha t}E\big[y'(t)P_{\alpha}y(t)
-y'(t)\int_0^{h_r}\Pi_{\alpha}(\theta)\hat{y}(t|t+\theta)d\theta\nonumber\\
&&-\int_0^{h_r}\tilde{y}'(t|t+\theta)\Pi_{\alpha}(\theta)
\tilde{y}(t|t+\theta)d\theta\big]\nonumber\\
&=&e^{\alpha t}E\big[y'(t)\hat{P}_{\alpha}y(t)\big]\geq 0,\label{zz8}
\end{eqnarray}
where $\tilde{y}(t|t+\theta)=y(t)-\hat{y}(t|t+\theta),$ and $\hat{P}_{\alpha}>0$ is used in the last equality. It is clear that $~V(t,y(t))\rightarrow \infty$ if $E\|y(t)\|^2\rightarrow \infty$ from (\ref{zz8}). By taking time derivative along the dynamic of the stochastic system (\ref{i3}) and combining with (\ref{x6})-(\ref{x10}), we have
\begin{eqnarray}
&&\dot{V}(t,y(t))\nonumber\\
&=&e^{\alpha t}E\Big\{y'(t)\big[A'P_{\alpha}+P_{\alpha}
A+\alpha P_{\alpha}-\Pi_{\alpha}(h_r)\big]y(t)\nonumber\\
&&+u(t)'\sum_{i=0}^{r}\bar{B}_i'e^{-A'h_i}P_{\alpha}(t)e^{-Ah_i}\bar{B}_i u(t)\nonumber\\
&&-y'(t)\int_t^{t+h_r}\Big[\frac{\partial}{\partial
t}\Pi_{\alpha}(\theta-t)+\Pi_{\alpha}(\theta-t)A
+A'\Pi_{\alpha}(\theta-t)\nonumber\\
&&+\alpha\Pi_{\alpha}(\theta-t)\Big]\hat{y}(t|\theta)d\theta
-u'(t)\Omega_{\alpha} K_{\alpha} y(t)\nonumber\\
&&-y'(t) K_{\alpha}'\Omega_{\alpha}'u(t)
+y'(t)\Pi_{\alpha}(0)\hat{y}(t|t) \Big\}\nonumber\\
&=&-e^{\alpha t}E\Big[y'(t)y(t)+u(t)' u(t)\Big]\leq
0.\label{zz14}
\end{eqnarray}
Thus from (\ref{zz14}), we know $V(t,y(t))$ is nonincreasing, and thus $V(t,y(t))\leq V(0,y(0))$. Therefore,
$\lim_{t\rightarrow\infty}V(t,y(t))$ exists.

Integrating on both sides of (\ref{zz14}) from $t$ to $t+T$ yields
\begin{eqnarray}
&&V(t+T,y(t+T))-V(t,y(t))\nonumber\\
&=&-\int_{t}^{t+T}e^{\alpha s}E\big[y'(s)y(s)+u'(s)u(s)\big]ds\nonumber\\
&=&-\int_{t}^{t+T}e^{\alpha s}E\big[y'(s)y(s)
+\hat{y}'(s|s)K_{\alpha}'K_{\alpha}\hat{y}(s|s)\big]ds.\nonumber
\end{eqnarray}

Now we consider the following cost function,
\begin{eqnarray}
E\int_{t}^{t+T}e^{\alpha s}\Big[{y}'(s)y(s)+{u}'(s)u(s)\Big]ds. \label{lele1}
\end{eqnarray}
By applying Theorem \ref{theorem}, the optimal controller to minimize (\ref{lele1}) subjected to system (\ref{i3}) is given as
$u^*(s)=K_{\alpha}(s)\hat{y}^*(s|s)$, where $K_{\alpha}(s)$ is given by (\ref{d4})-(\ref{d5}) with $Q=I,R=I.$
$y^*(s)$ is the corresponding state trajectory. Accordingly, the optimal cost of (\ref{lele1}) is given by
\begin{eqnarray}
&&E\int_{t}^{t+T}e^{\alpha s}\big[{y^*}'(s)y^*(s)+{u^*}'(s)u^*(s)\big]ds\nonumber \\
&=&e^{\alpha t}E\big[y'(t)P_{\alpha}(0)y(t)-y'(t)\int_0^{h_r}\Pi_{\alpha}(0,\theta)\hat{y}(t|t+\theta)d\theta\big],\nonumber
\end{eqnarray}
Therefore, we have
\begin{eqnarray}
&&V(t+T,y(t+T))-V(t,y(t))\nonumber\\
&=&-E\int_{t}^{t+T}e^{\alpha s}\big[y'(s)y(s)+\hat{y}'(s|s)K_{\alpha}'K_{\alpha} \hat{y}(s|s)\big]ds\nonumber \\
&\leq&-E\int_{t}^{t+T}e^{\alpha s}\big[{y^*}'(s)y^*(s)+{u^*}'(s)u^*(s)\big]ds\nonumber\\
&=&-e^{\alpha t}E\big[y'(t)P_{\alpha}(0)y(t)
-y'(t)\int_0^{h_r}\Pi_{\alpha}(0,\theta)\hat{y}(t|t+\theta)d\theta\big]\nonumber\\
&\leq&0.\label{kkV}
\end{eqnarray}

Note
\begin{eqnarray}
&&\lim_{t\rightarrow\infty}[V(t+T,y(t+T))-V(t,y(t))]\nonumber\\
&=&\lim_{t\rightarrow\infty}V(t+T,y(t+T))-\lim_{t\rightarrow\infty}V(t,y(t))=0,\nonumber
\end{eqnarray}
it follows from (\ref{kkV}) that
\begin{eqnarray}
0&=&\lim_{t\rightarrow\infty}e^{\alpha t}E\Big[y'(t)P_{\alpha}(0)y(t)\nonumber\\
&&-y'(t)\int_0^{h_r}\Pi_{\alpha}(0,\theta)\hat{y}(t|t+\theta)d\theta\Big].\nonumber
\end{eqnarray}
Further, since
\begin{eqnarray}0&\leq&
E\big[y'(t)\hat{P}_{\alpha}(0)y(t)\big]\nonumber\\
&\leq&
E\big[y'(t)P_{\alpha}(0)y(t)-y'(t)\int_0^{h_r}\Pi_{\alpha}(0,\theta)\hat{y}(t|t+\theta)d\theta\big],\nonumber
\end{eqnarray}
it follows that
\begin{eqnarray}
\lim_{t\rightarrow\infty}e^{\alpha t}E\big[y'(t)\hat{P}_{\alpha}(0)y(t)\big]=0.\label{zz16}
\end{eqnarray}
Now we are in the position to show that $\hat{P}_{\alpha}(0)>0.$ If this
is not the case, there would exist $z\neq0,$ such that
$z'\hat{P}_{\alpha}(0)z=0$. Consider the closed-loop system
$dy(t)=[Ay(t)+BK_{\alpha}\hat{y}(t|t)]dt+\sum_{i=0}^{r}e^{-Ah_i}\bar{B}_iK_{\alpha}\hat{y}(t|t)dw_i(t+h_i)$ with initial value $y(0)=z$.
Return to (\ref{zz18}), one has
\begin{eqnarray}
E\int_{0}^{T}[{y^*}'(t)y^*(t)+{u^*}'(t)u^*(t)]dt=z'\hat{P}(0)z=0.\nonumber
\end{eqnarray}
Together with Assumption 1, one has
\begin{eqnarray}
y^*(t)=0,u^*(t)=0,t\geq 0, a.s._.\nonumber
\end{eqnarray}
The system (\ref{i4}) is thus reduced to
\begin{eqnarray}dy^*(t)=Ay^*(t)dt,
y^*(t)=0, a.s._, ~t\in[0,T].\nonumber
\end{eqnarray}
In view of the observability of system $(A,I)$, it
yields that $z=0,$ which is a contradiction. Thus,
$\hat{P}_{\alpha}(0)>0.$ Together with (\ref{zz16}),
we have
$$\lim_{t\rightarrow\infty}e^{\alpha t}E\|y(t)\|^2=0.$$

Using the fact that $u(t)=K_{\alpha}\hat{y}(t|t),$ it is immediately obtained that $$\lim_{t\rightarrow\infty}e^{\alpha t}E\|u(t)\|^2=0.$$
Thus $\lim_{t\rightarrow\infty}e^{\alpha t}E\|x(t)\|^2=0$ follows from (\ref{i3}).
The exponential mean-square
stability of system (\ref{i1}) follows.
The proof is now completed.


\end{document}